\documentclass[reqno, a4paper, 11pt]{amsart}
\usepackage{amssymb, amsthm, amsmath}

\numberwithin{equation}{section}

\renewcommand{\d}{\mathrm{d}}
\newcommand{\x}{\ma{x}}
\newcommand{\w}{\ma{w}}
\newcommand{\y}{\ma{y}}

\renewcommand{\v}{\ma{v}}
\renewcommand{\u}{\ma{u}}
\renewcommand{\t}{\ma{t}}

\newcommand{\mla}{\boldsymbol{\lambda}}
\newcommand{\mbeta}{\boldsymbol{\beta}}

\newcommand{\TT}{\mathcal{T}}
\newcommand{\SSS}{\mathcal{S}}
\newcommand{\WW}{\mathcal{W}}
\newcommand{\R}{\mathbb{R}}

\newcommand{\Z}{\mathbb{Z}}
\newcommand{\N}{\mathbb{N}}
\newcommand{\Q}{\mathbb{Q}}

\newcommand{\bfP}{\mathbb{P}}
\newcommand{\A}{\mathbb{A}}

\newcommand{\ov}{\overline}
\newcommand{\ma}{\mathbf}

\newcommand{\mcal}{\mathcal}
\renewcommand{\le}{\leqslant}
\renewcommand{\ge}{\geqslant}
\renewcommand{\leq}{\leqslant}
\renewcommand{\geq}{\geqslant}
\newcommand{\ben}{\begin{enumerate}}
\newcommand{\een}{\end{enumerate}}
\newcommand{\eit}{\begin{itemize}}
\newcommand{\beq}{\begin{equation}}
\newcommand{\eeq}{\end{equation}}
\newcommand{\ve}{\varepsilon}

\newcommand{\al}{\alpha}
\newcommand{\D}{\Delta}

\newcommand{\la}{\lambda}
\newcommand{\sfl}{\mathsf{\Lambda}}

\renewcommand{\mod}{\hspace{-0.25cm}\pmod}

\newcommand{\lab}{\label}

\newcommand{\hcf}{\mathrm{gcd}}
\renewcommand{\mod}{\hspace{-1mm}\pmod}
\newcommand{\colt}[2]{\genfrac{}{}{0pt}{1}{#1}{#2}}
\newcommand{\tstack}[3]{\colt{#1}{\colt{#2}{#3}}}

\renewcommand{\b}[1]{{\bf #1}}

\newcommand{\ep}{\varepsilon}

\newtheorem*{thm*}{Theorem}
\newtheorem{lem}{Lemma}
\newtheorem{pro}{Proposition}
\newtheorem*{pro*}{Proposition}

\newtheorem*{cor*}{Corollary}

\newtheorem*{con}{Conjecture}
\theoremstyle{definition}

\newtheorem*{notat}{Notation}

\DeclareMathOperator{\sing}{sing}
\DeclareMathOperator{\meas}{meas}

\begin{document}

\title{Counting rational points on cubic hypersurfaces}

\author{T.D. Browning}
\address{School of Mathematics,  University of Bristol, Bristol BS8 1TW}
\email{t.d.browning@bristol.ac.uk}

\date{\today}

\begin{abstract}
Let $X\subset \bfP^{N}$ be a geometrically integral cubic hypersurface
defined over $\Q$, with singular locus of dimension $\leq \dim X-4$. 
Then the main result in this paper is a proof of the fact that $X(\Q)$
contains $O_{\ve,X}(B^{\dim X+\ve})$ points of height at most $B$.
\end{abstract}

\maketitle

\section{Introduction}

Let $C\in\Z[x_1,\ldots,x_n]$ be an absolutely irreducible
cubic form, defining a hypersurface $X_C \subset \bfP^{n-1}$.
The primary goal of this paper is to
investigate the density of rational points on $X_C$.  Given any
rational point  $x=[\x] \in \bfP^{n-1}(\Q)$, with 
$\x=(x_1,\ldots,x_n) \in \Z^{n}$ and $\hcf(x_1,\ldots,x_n)=1$, we
write  $H(x):=\max_{1\leq i \leq n}|x_i|$. Let 
$$
N_{X_F}(P):=\#\{x \in X_F\cap\bfP^{n-1}(\Q): ~H(x) \leq P\},
$$
for any $P\geq 1$ and any hypersurface $X_F\subset \bfP^{n-1}$ defined 
by a form $F\in\Z[\x]$ in $n$ variables.  This paper is motivated by the following basic
conjecture,  due to Heath-Brown \cite[Conjecture 2]{annal}.

\begin{con}  Let $\ve>0$ and suppose that 
$F \in \Z[x_1,\ldots,x_n]$ is an absolutely irreducible form of degree
$d\geq 2$. Then we have
$$
N_{X_F}(P)=O_{d,\ve,n}(P^{n-2+\ve}).
$$
\end{con}

We will henceforth refer to this conjecture as the ``uniform dimension
growth conjecture''. There is a version of the conjecture in which one
allows an arbitrary dependence on the coefficients of the form under
consideration. We will refer to this as the 
``dimension growth conjecture''. These conjectures are essentially
best possible, as examples of the shape 
$$
F(\x)=x_1^d+x_2(x_3^{d-1}+\cdots+x_n^{d-1})
$$
show. Here, one obtains $N_{X_F}(P)\gg P^{n-2}$ by considering
rational points of the shape $[0,0,x_3,\ldots,x_n]$.  The dimension
growth conjectures have received
considerable attention in recent years, to the extent that the uniform
dimension growth conjecture is known to hold 
when $F$ is non-singular, when $d=2$ or $n\leq 5$, and when $d\geq 4$.
This is one of the major outcomes of the body of work \cite{n-2, smoothII,
  annal, bonn}. Thus the single outstanding case concerns singular absolutely
irreducible cubic forms $C\in\Z[x_1,\ldots,x_n]$, with $n\geq 6$.
In this setting, the best general result available 
is due to Salberger \cite{bonn}, who has shown that
\begin{equation}
  \label{eq:record}
N_{X_C}(P) \ll_{\ve,n}P^{n-3+2/\sqrt{3}+\ve},  
\end{equation}
for any $\ve>0$. Note that $-1+2/\sqrt{3}=0.1547\ldots$.

Let $s(C)\in\Z\cap[-1,n-2]$ denote the projective dimension of the
singular locus of the hypersurface $X_C$. 
When $s(C)=-1$ we have already seen that the uniform dimension growth
conjecture holds. One might therefore hope to improve on \eqref{eq:record} when $s(C)$
is not too large compared with $n$. 
This is the point of view adopted by Salberger \cite{corput}, who has
recently established the validity of the conjecture for $n\geq
8+s(C)$. Salberger's result is based on refining an argument involving
the $q$-analogue of van der Corput's method for exponential sums,
developed by Heath-Brown \cite{hb-india}.  As observed independently
by the author and Salberger, the weaker bound $n\geq 11+s(C)$ is
readily derived from Heath-Brown's work.
At the expense of uniformity in the implied constant we will be able
to improve on all of these results. Let $\|C\|$ denote the maximum modulus  of
the coefficients of the underlying cubic form $C$. Then the following is
our main result.

\begin{thm*}
Let $X_C\subset \bfP^{n-1}$ be a geometrically integral cubic
hypersurface, and let $\ve>0$.
Assume that 
$
n\geq 6+s(C).
$
Then there exists a positive number $\theta=O_n(1)$ such that
$$
N_{X_C}(P)\ll_{\ve,n}\|C\|^{\theta} P^{n-2+\ve}.
$$
In particular, the dimension growth conjecture holds for $n\geq 6+s(C).$
\end{thm*}

With more work, the value of $\theta$ could be made explicit in the
statement of the theorem. We have not attempted to do so here,
however, being content to show that there is at worst
polynomial dependence on $\|C\|$.  An inspection of the proof of
the theorem  reveals that we obtain a small improvement on the 
exponent $n-2$ as soon as $n>6+s(C)$.

We will actually establish a version of the theorem for the
number of integral points in certain expanding 
regions that lie on arbitrary affine cubic hypersurfaces $g=0$.
The proof of this estimate will be established by
induction on the dimension $s(g_0)$ of the singular locus of the
projective hypersurface $g_0=0$, where $g_0$ denotes the
cubic part of the polynomial $g$.  
The idea is to use hyperplane sections to reduce 
consideration to a family of hypersurfaces in $\A_{\Q}^{n-1}$, for which the dimension of
the relevant singular locus is reduced by $1$.
The main work comes in having to handle the inductive base
$s(g_0)=-1$, and this will be dealt with by an application of the
Hardy--Littlewood circle method.

\begin{notat}
Throughout our work $\N$ will denote the set of positive
integers.  For any $\al\in \R$, we will follow common convention and
write $e(\al):=e^{2\pi i\al}$ and $e_q(\al):=e^{2\pi i\al/q}$.  We
will allow both the small positive quantity $\ve$ and the constant
$\theta$ to vary from time to
time, so that we may write $P^{\ve}\log P\ll_{\ve}P^{\ve}$ and
$H^{2\theta}\leq H^{\theta}$, for
example.  All of the implied
constants in our work are allowed to depend implicitly on $\ve$ and
$n$, and $\theta$ will 
always be bounded in terms of $n$ alone. Finally, we will write
$|\ma{z}|:=\max_{1\leq i\leq n}|z_i|$ for the norm of any vector
$\ma{z}\in\R^{n}$, and we will use the notation $A \asymp B$ to denote that 
$A\ll B \ll A$. 
\end{notat}

\section{Preliminaries}\label{s:act}

It will be convenient to work with infinitely differentiable weight functions
$w: \R^n \rightarrow  \R_{\geq 0}$, with compact support. Given such a function $w$, we set
$R(w)$ to be the smallest $R$
such that $w$ is supported in the hypercube $[-R,R]^n$, and we let
$$
R_j(w):=\max\Big\{ \Big|
\frac{\partial^{j_1+\cdots+j_n}w(\x)}{\partial^{j_1}x_1\cdots
\partial^{j_n}x_n}\Big|: ~\x \in \R^n, ~j_1+\cdots+j_n=j\Big\},
$$
for each integer $j \geq 0$.  Let constants $c_n$ and $c_{n,j}$ be
given, and define $\WW_n$ to be 
the set of infinitely differentiable functions $w: \R^n \rightarrow
\R_{\geq 0}$ of compact support, such that $R(w)\leq c_n$ and
$R_j(w)\leq c_{n,j}$ for all $j\ge 0$.  In future all our order
constants will be allowed to depend on $c_n$ and the $c_{n,j}$,
without further comment.

Let $g\in \Z[x_1,\ldots,x_n]$ be a cubic polynomial and let $w\in\WW_n$.
All our efforts are centred upon determining  the asymptotic
behaviour of the quantity 
$$
N_{w}(g;P):= \sum_{\colt{\x\in \Z^n}{g(\x)=0}} w(\x/P),  
$$
as $P\rightarrow \infty$. 
For this we will employ the form of the Hardy--Littlewood circle
method developed by Heath-Brown \cite{hb-10}, which incorporates a
single Kloosterman refinement. Define the cubic exponential sum
\begin{equation}\label{eq:Tcubic}
\TT(\al)=\TT_P(\al;g,w):=\sum_{\x\in\Z^n}w(\x/P)e(\al g(\x)),
\end{equation}
for $P\geq 1$. Then $\TT(\al)$ converges absolutely,
and for any $Q\geq 1$ we have 
$$
N_{w}(g;P)=\int_0^1 \TT(\al)\d\al=
\int_{\frac{-1}{1+Q}}^{1-\frac{1}{1+Q}} \TT(\al)\d\al.
$$
In \cite{hb-10} Heath-Brown proceeds to break the interval 
$[\frac{-1}{1+Q}, 1-\frac{1}{1+Q}]$ 
according to the Farey dissection of order $Q$. This ultimately yields
\begin{equation}
  \label{eq:Nw'}
\begin{split}
N_{w}(g;P)=&\sum_{q\leq Q}\int_{\frac{-1}{qQ}}^{\frac{1}{qQ}} \SSS_0(q;z)\d
 z +O\big(Q^{-2}E_{w}(g;P,Q)\big),
\end{split}
\end{equation}
for any $Q\geq 1$, where 
$$
  E_{w}(g;P,Q):=
\sum_{q\leq Q}\sum_{|u|\leq
\frac{q}{2}}\frac{\max_{\frac{1}{2}\leq qQ|z|\leq 1}|\SSS_u(q;z)|}{1+|u|},
$$
and 
\begin{equation}
  \label{eq:Su}
\SSS_u(q;z):=\sum_{\colt{a=1}{\hcf(a,q)=1}}^q e_q(\ov{a}u)\TT(a/q+z).
\end{equation}
This is \cite[Lemma 7]{hb-10}. 
We will find that our work is optimised by taking $Q=P^{3/2}$ in
\eqref{eq:Nw'}.  

Recall the notation $\|g\|$ for the maximum modulus of the
coefficients of $g$.  We will follow the convention that $g_0$ denotes the 
homogeneous cubic part of $g$. Given $P \geq 1$, we will need to work with the function
\begin{equation}  \label{eq:norm_P}
\|g\|_P:= \|P^{-3}g(Px_1, \ldots, Px_n)\|.
\end{equation}
It is clear that
$$
\|g_0\|\leq \|g\|_P \leq \|g\|,
$$
with equality if $g$ is homogeneous.
The bulk of this paper will be spent establishing the
following result, from which the proof of the theorem will
flow rather swiftly.

\begin{pro}\label{main'}
Suppose we are given $w\in \WW_n$, $\ve>0$ and $H\ge 1$.
Let $g\in  \Z[x_1,\ldots,x_n]$ be a cubic polynomial for which $g_0$ is non-singular.
Assume that $n\geq 5$.
Then there exists a positive number $\theta$ such that
\begin{equation}
  \label{eq:gen1}
N_{w}(g;P)\ll H^{\theta}P^{n-2+\ve},
\end{equation}
for any $H \geq \|g\|_P.$
\end{pro}

The proof of this result relies upon \eqref{eq:Nw'} in a crucial way,
and will involve two basic estimates
for $\SSS_u(q;z)$. The first uses repeated Weyl
differencing, and is based on the approach taken by Davenport in
\cite{D16}. This will be the subject of \S \ref{s:cubic_weyl}.
The second estimate is based on an application of the Poisson
summation formula, and in particular, 
the treatment of cubic exponential sums appearing in the author's recent joint work with Heath-Brown \cite{41}.
This will be the focus of \S\ref{s:cubic}. 
Finally, in \S \ref{s:finale} we will stitch all of 
this together in order to complete the proof of Proposition~\ref{main'}.

We end this section by showing how the theorem follows from
Proposition~\ref{main'}.   
Let $g\in  \Z[x_1,\ldots,x_n]$ be a cubic polynomial with cubic part $g_0$.
Let $w\in \WW_n$ and $\ve>0$ be given. We will show that
there exists $\theta>0$ such that
\eqref{eq:gen1} holds, provided that $H \geq \|g\|_P$ and 
\begin{equation}
  \label{eq:gen2}
n\geq 6+ s(g_0).
\end{equation}
Before establishing this claim, let us indicate how this suffices for
the statement of the theorem. Define 
$$
w_1(\x):=\prod_{i=1}^n \gamma(x_i),
$$
where $\gamma:\R\rightarrow \R_{\geq 0}$ is given by 
$$
\gamma(x)=\left\{
\begin{array}{ll}
e^{-1/(1-x^2)}, & \mbox{if $|x|<1$},\\
0, & \mbox{if $|x|\geq 1$}.
\end{array}
\right.
$$
It is clear that $w_1 \in \WW_n$.
Let $C\subset \Z[\x]$ be a cubic form defining a cubic hypersurface
$X_C\subset \bfP^{n-1}$. Then it follows from the above claim that 
$$
N_{X_C}(P)\ll N_{w_1}(C;P) \ll H^{\theta}P^{n-2+\ve},
$$
provided that $n\geq 6+ s(C_0)=6+s(C)$ and $H\geq \|C\|_P=\|C\|$. 
This therefore establishes the theorem subject to the claim.

To confirm the claim we will argue by induction on $s(g_0)$, the
base case $s(g_0)=-1$ being taken care of by Proposition \ref{main'}. 
To handle the inductive step we will use a
simpler version of the argument based on hyperplane sections developed
in \cite[\S 4]{41}. 
Let $w\in \WW_n$ and let $g \in \Z[x_1,\ldots,x_n]$ be a cubic
polynomial such that $s(g_0)\geq 0$.  
Let $P\geq 1$ and let $H$ be such that $\|g\|_P\leq H.$ 
Our plan will be to use hyperplane sections, in order to reduce the problem to a
consideration of cubic polynomials in only $n-1$ variables,
whose cubic part defines a hypersurface with  singular locus of dimension $s(g_0)-1$.
According to \cite[Lemma 5]{41} in the special case $d=3$ and $r=0$, there
exists a primitive vector  $\ma{m}\in\Z^n$, with $|\ma{m}|\ll 1$,  
such that 
\begin{equation}
  \label{eq:pizza}
\dim \sing (X_{g_0}\cap H_{\ma{m}})= s(g_0)-1.
\end{equation}
Here $X_{g_0}$ is the hypersurface defined by $g_0=0$ and $H_{\ma{m}}$
is the hyperplane $\ma{m}.\x=0$.
In order to apply the induction hypothesis we will sum over affine
hyperplane sections $\ma{m}.\x=k$, for integers $k\ll P$.
This gives 
\begin{equation}
\lab{eq:ind}
N_{w}(g;P)
= \sum_{k\ll P} 
\sum_{\tstack{\x\in\Z^n}{\ma{m}.\x=k}{g(\x)=0}} w(\x/P) \\
=\sum_{k\ll P} \mcal{N}_k,
\end{equation}
say. 
Now $\mcal{N}_k$ is zero unless there exists a vector $\ma{t}\in \Z^n$
such that $\ma{m}.\ma{t}=k$ and $|\ma{t}|\ll P$. Let us fix such a
choice of vector, and write $\x=\ma{t}+\y$ in $\mcal{N}_k$.
Then clearly $\ma{m}.\x=k$ if and only if $\ma{m}.\y=0$. 
This condition defines a lattice $\sfl\subseteq \Z^n$ of rank $n-1$ and
determinant $|\ma{m}|$, by part (i) of \cite[Lemma~1]{annal}. We now choose a
basis $\ma{e}_1,\ldots, \ma{e}_{n-1}$ for $\sfl$. Then
$1\leq |\ma{e}_i| \ll 1,$ for $1\leq i\leq n-1$.
Moreover, any of the vectors $\y$ we are interested in can be written as
$\y=\sum_{i=1}^{n-1}\la_i \ma{e}_i$, for 
$\mla=(\la_1,\ldots,\la_{n-1})\in\Z^{n-1}$ such that $\mla \ll  P$.

Putting all of this together, we conclude that 
\begin{align}
  \label{eq:ind'}
  \mcal{N}_k&= \sum_{\colt{\mla\in \Z^{n-1}}{h(\mla)=0}} w_0\big(\mla /P\big)
=N_{w_0}(h;P),
\end{align}
where 
$$
h(\u):=g\Big(\t+\sum_{i=1}^{n-1}u_i \ma{e}_i\Big), \quad
w_0(\ma{u}):=
w\Big(P^{-1}\t+\sum_{i=1}^{n-1}u_i \ma{e}_i\Big),
$$
and $\ma{u}=(u_1,\ldots,u_{n-1})$. Note that $h$ is a cubic polynomial
in only $n-1$ variables. We need to show that the induction hypothesis
can be applied to estimate  $N_{w_0}(h;P)$.
Now it is trivial to see that $w_0\in\WW_{n-1}$ for our choice of $\ma{m}$ and
$\ma{t}$, and furthermore, that 
$\|h\|_{P} \ll H $. Finally, the argument in \cite[\S 4]{41}
ensures that $s(h_0)=  s(g_0)-1$, under the assumption that
\eqref{eq:pizza} holds. 
On applying the induction hypothesis in \eqref{eq:ind'}, and combining
it with \eqref{eq:ind}, we therefore deduce that 
\begin{align*}
N_w(g;P) \ll H^{\theta} \sum_{k\ll P}
P^{n-3+\ve} \ll H^{\theta} P^{n-2+\ve}.
\end{align*}
This completes the proof of \eqref{eq:gen1} subject to  \eqref{eq:gen2}.

\section{Estimating $\SSS_u(q;z)$: Weyl differencing}\label{s:cubic_weyl}

In this section we will establish an estimate for the 
exponential sum \eqref{eq:Su} by arguing along the lines of 
Davenport \cite{D16}. The fact that we are working with possibly 
non-homogeneous polynomials makes no difference to the opening steps
of the argument. It will be convenient to 
draw upon Heath-Brown's recent reworking of
Davenport's approach \cite{14}, where possible.
Throughout this section we will suppose that $s(g_0)=-1$, so that the cubic part
$g_0$ of $g$ is non-singular.

We will sum trivially over the numerator in \eqref{eq:Su}, giving 
\begin{equation}
  \label{eq:Su_triv}
  |\SSS_u(q;z)|\leq q \max_{\colt{1\leq a \leq q}{\hcf(a,q)=1}}|\TT(a/q+z)|.
\end{equation}
Our interest lies with values of $q\leq Q$ and $|z|\leq (qQ)^{-1}$.
As indicated in \S \ref{s:act}, the final analysis will be
optimised by taking $Q=P^{3/2}$. In particular, we may henceforth
assume that 
\begin{equation}
  \label{eq:assume_z}
|z|\leq q^{-1}P^{-3/2}.  
\end{equation}
The purpose of this section is to establish the following result.

\begin{pro}\lab{pro:Sbirch}
Let $q\in \N$ such that $q\leq P^{3/2}$, and let $z\in \R$ such that 
\eqref{eq:assume_z} holds.  
Let $w\in\WW_n$, let $\ve>0$ and let $g\in \Z[x_1,\ldots,x_n]$ be a
cubic polynomial such that $g_0$ is non-singular.  Then we have
$$
\SSS_u(q;z)
\ll \|g_0\|^{n/8} q^{1-n/8} P^{n+\ep}
\min\big\{1, (|z|P^{3})^{-n/8}\big\}.
$$
\end{pro}

Recall the definition \eqref{eq:Tcubic} of $\TT(\al)=\TT_P(\al;g,w)$,
for $w \in \WW_n$. The central idea in Davenport's approach is an application of Weyl 
differencing. The first step in this process produces the bound
$$
|\TT(\alpha)|^2\ll \sum_{\w\ll P}\Big|\sum_{\x}
w\big((\x+\w)/P\big)w(\x/P)
e\big(\al (g(\x+\w)-g(\x))\big)\Big|.
$$
An application of Cauchy's inequality now yields
\begin{equation}\label{22-weyl3}
|\TT(\alpha)|^4\ll P^n\sum_{\w,\x\ll P}\Big|\sum_{\y\in \Z^n}
w_{\w,\x}(\y)
e\big(\alpha G(\w,\x;\y)\big)\Big|,
\end{equation}
where 
\begin{equation}
  \label{eq:G()}
  G(\w,\x;\y):=g(\w+\x+\y)-g(\w+\y)-g(\x+\y)+g(\y),
\end{equation}
and 
$$
w_{\w,\x}(\y)=w\big((\w+\x+\y)/P\big)
w\big((\w+\y)/P\big)w\big((\x+\y)/P\big)w(\y/P).
$$

Recall our notation $g_0$ for the homogeneous cubic part of 
the cubic polynomial $g$. Suppose that
$$
g_0(x_1,\ldots,x_n)=\sum_{i,j,k=1}^n c_{ijk}x_ix_jx_k,
$$
in which the coefficients $c_{ijk}$ are symmetric in the indices
$i,j,k$. On replacing $g$ by $6g$,  we may assume that the $c_{ijk}$
are all integral. If we now define the bilinear forms
$$
B_i(\w;\x):=\sum_{j,k=1}^n c_{ijk}w_jx_k,
$$
for $1\leq i\leq n$, then we find that
$$
G(\w,\x;\y)=6\sum_{i=1}^n y_iB_i(\w;\x)+\Gamma(\w,\x),
$$
where $G$ is given by \eqref{eq:G()} and $\Gamma(\w,\x)$ is 
independent of $\y$.
It therefore follows from \eqref{22-weyl3} that
\begin{align*}
|\TT(\alpha)|^4
&\ll 
P^n\sum_{\w,\x\ll P}\Big|\sum_{\y\in \Z^n}
w_{\w,\x}(\y)
e\Big(6\alpha \sum_{i=1}^n y_i B_i(\w;\x)\Big)\Big|.
\end{align*}
On combining the standard estimate for linear exponential sums with
partial summation, we deduce that
\begin{align*}
|\TT(\alpha)|^4
&\ll P^{n} \sum_{\w,\x\ll P}\,
\prod_{i=1}^n
\min\{P, \| 6\alpha B_i(\w;\x)\|^{-1}\}.
\end{align*}
The aim is now to establish a link between this bound 
and the density of integer solutions to the 
system of simultaneous bilinear equations
\begin{equation}\label{bileq}
B_i(\x;\b{y})=0,\quad (1\leq i\leq n).
\end{equation}
This is described in detail by Heath-Brown
\cite[\S 2]{14}. Following this more or less verbatim we obtain
$$
|\TT(a/q+z)|^4\ll \frac{P^{2n}(\log P)^n}{Z^{2n}} 
\#\Big\{(\w,\x)\in\Z^{2n}:
\begin{array}{l}
\w,\x\ll ZP,\\ 
\mbox{\eqref{bileq} holds}
\end{array}
\Big\},
$$ 
for any $Z\in \R$ such that 
$$
0<Z<1,\quad Z^2<(12cq|z|P^2)^{-1},\quad Z^2<P/(2q),
$$
and
$$
Z^2<\max\Big\{\frac{q}{6cP^2}\,,\,qP|z|\Big\}.
$$
Here, $c=\sum |c_{ijk}|$, where $c_{ijk}$ are the coefficients of
$g_0$. In particular, we clearly have $\|g_0\|\ll c \ll \|g_0\|$.

Now it is not hard to see that the system of equations 
\eqref{bileq} is just the system $\ma{H}_{g_0}(\x)\y=\ma{0}$,
where 
$$
\ma{H}_F(\x):=\Big\{\frac{\partial^2 F}{\partial x_i
\partial x_j}\Big\}_{1\leq i,j\leq n}  
$$
is the Hessian matrix formed from the
second order partial derivatives of any form $F\in \Z[\x]$. 
It follows from \cite[Lemma 1]{hb-10} that the variety cut out by \eqref{bileq} has
dimension $n$ in $\A^{2n}$.
An application of \cite[Eq. (2.3)]{n-2} now yields
\begin{align*}
|\TT(\al)|^4
&\ll \frac{P^{2n}(\log P)^n}{Z^{2n}} 
(ZP)^{n}\ll
Z^{-n}P^{3n} (\log P)^n,
\end{align*}
provided that $Z \geq P^{-1}$.  This bound clearly holds trivially
when $Z<P^{-1}$.
We will need to choose $Z$ as large as possible, given the constraints
above.  The choice 
$$
Z=\frac{1}{2}\min\Big\{1\,,\, \frac{1}{12cq|z|P^2}\,,\,
\frac{P}{2q}\,,\,\max\big\{\frac{q}{6cP^2},
q|z|P\big\}\Big\}^{1/2},
$$
is clearly satisfactory. On taking this value, we therefore deduce that
\begin{align*}
|\TT(\al)|^4 & \ll_{\ve} \|g_0\|^{n/2}P^{3n+\ep}
\Big(1+q|z|P^2+qP^{-1}
+\frac{1}{q}\min\big\{P^2,\frac{1}{|z|P}\big\}\Big)^{n/2}\\
& = \|g_0\|^{n/2} P^{4n+\ep}
\Big(P^{-2}+q|z|+qP^{-3}
+\frac{1}{q}\min\big\{1,\frac{1}{|z|P^3}\big\}\Big)^{n/2},
\end{align*}
for any $\ve>0$. We complete the proof of Proposition
\ref{pro:Sbirch} by substituting this into \eqref{eq:Su_triv}.

\section{Estimating $\SSS_u(q;z)$: Poisson summation}\label{s:cubic}

The goal of this section is to give an alternative treatment of the
cubic exponential sum \eqref{eq:Su}, for cubic
polynomials $g\in \Z[x_1,\ldots, x_n]$ such that the cubic part $g_0$
is non-singular. Our treatment is based on drawing together
ideas already present in the author's joint work with Heath-Brown
\cite{41}, and the Kloosterman refinement carried out by Heath-Brown
\cite{hb-10} in the context of non-singular cubic forms.

We will continue to assume that
$\al=a/q+z$, with $q\leq P^{3/2}$ and $z$ satisfying 
\eqref{eq:assume_z}.  We will write $q=bc^2d$, where
\begin{equation}
  \label{eq:bcd}
b_1:=\prod_{p \| q}p, 
\quad
b_2:=\prod_{p^2\| q}p, 
\quad
d:=\prod_{\colt{p^e\| q}{e\geq 3, ~2\nmid e}}p
\end{equation}
and $b=b_1b_2^2$.  It is not hard to see that $d$ divides $c$, and that 
there exist a divisor $d_0$ of $d$ such that
$d_0^{-1}d^{-1}c$ is a square-full integer. Moreover, 
$\hcf(b,c^2d)=1$.  Finally, we recall the definition of the
function \eqref{eq:norm_P} for $P \geq 1$.
The following is our main estimate for $\SSS_u(q,z)$.

\begin{pro}\label{pro:S}
Let $w\in\WW_n$, let $\ve>0$ and let $g\in \Z[x_1,\ldots,x_n]$ be a
cubic polynomial such that $g_0$ is non-singular
and $\|g\|_P\leq H$, for some
$H\leq P.$ Let $q=bc^2d$, in the notation of \eqref{eq:bcd}. Define
\begin{equation}
  \label{eq:V}
V:=qP^{-1} \max\{1,\sqrt{|z|P^3}\},
\end{equation}
and
\begin{equation}
  \label{eq:W}
  W:=V+ (c^2d)^{1/3}.
\end{equation}
Then there exists a positive number $\theta$ such that
$$
\SSS_u(q;z)\ll H^\theta 
q^{-n/2+1} P^{n+\ve} \big(W^n M_1+\min\{W^n,M_2,M_3\}\big),
$$
where
\begin{align*}
M_1
&:=\max_{0< N \ll
  (HP)^{\theta}}\frac{\hcf(b_1,N)^{1/2}}{b_1^{1/2}},
\end{align*}
and 
$$
M_2:=c^{n}\Big(1+\frac{V}{c}\Big)^{n-3/2},\quad 
M_3:=V^n\Big(1+\frac{c^2d}{V^3}\Big)^{n/2}.
$$
\end{pro}

Note that by taking $M_1\leq 1$ and $\min\{W^n, M_2,M_3\}\leq W^n$ in this
estimate we retrieve what is given by summing trivially over
$a$ in \cite[Proposition~1]{41}.
Not surprisingly, our sharpening of this estimate is based
on taking advantage of cancellation in the summation over $a$.   

The first step in the proof of Proposition \ref{pro:S} involves introducing complete
exponential sums modulo $q$ via an application of Poisson summation. 
Thus it follows from summing over $a$ in \cite[Lemma 8]{41} that
$$
\SSS_u(q;z)=q^{-n}\sum_{\v\in\Z^n}S_u(q;\v)I(z;q^{-1}\v),
$$
where
\begin{equation}
  \label{eq:Sq}
S_u(q;\v):=\sum_{\colt{a=1}{\hcf(a,q)=1}}^q e_q(\ov{a}u)
\sum_{\y\bmod{q}}e_q(ag(\y)+\v.\y)
\end{equation}
and 
$$
  I(z;\mbeta):=\int w(\x/P)e(z g(\x)-\mbeta.\x)\d\x.
$$
Before discussing the complete exponential sums \eqref{eq:Sq}, we
first dispatch the integral $I(z;q^{-1}\v)$ that appears in our
formula for $\SSS_u(q;z)$.   For this we may apply \cite[Lemma
9]{41}. 
This leads to the conclusion that 
\begin{equation}  \label{eq:S2}
\begin{split}
  \SSS_u(q;z)&\ll P^{-N}+q^{-n}\int_{\x\ll P}
  \sum_{|\v-qz\nabla g(\x)|\leq P^\ve V}|S_u(q;\v)| \d\x\\
&\ll P^{-N}+q^{-n}P^n\max_{\v_0\ll H P}\sum_{|\v-\v_0|\leq P^\ve V}|S_u(q;\v)|.
\end{split}
\end{equation}
Here $N\in\N$ is arbitrary, $V$ is given by \eqref{eq:V} and we have
used the fact that 
$$
qz\nabla g(\x)\ll P^{-1}\|\nabla g(P\x)\| = P\|P^{-2}\nabla
g(P\x)\| \ll HP,
$$
for any $\x\ll P$ and $|z|\leq q^{-1}P^{-3/2}\leq (qP)^{-1}$.

The bound in \eqref{eq:S2} constitutes a key difference between our current
approach and that used in Heath-Brown's work on cubic forms
\cite{hb-10}. In the latter one uses instead the bound
$$
\SSS_u(q;z)\ll \frac{1}{P^{N}}+\frac{1}{q^{n}}\sum_{\v\ll HP^{1+\ve}}|S_u(q;\v)|
\meas\{\x\ll P:\,|\v_0(\x)-\v|\leq P^{\ve}V\},
$$
for $\v_0(\x)=qz\nabla g(\x)$.
In fact this difference is already capitalised upon in \cite{41}.  The point is that the approach in \cite{hb-10} 
requires information about the size of
the Hessian $\det\ma{H}_{g_0}(\x)$, and this is achieved by working
with weight functions that detect points lying very close to a
fixed point $\x_0 \in \R^n$ satisfying $g_0(\x_0)=0$, but 
which does not vanish on the
Hessian.  In our present investigation we seek an upper bound for 
integer solutions $\x \in \Z^n$ to the equation $g=0$, that have
modulus at most $P$, and not just those that lie sufficiently close to
$\x_0$. Thus we have found it convenient to adopt the device pioneered
in \cite{41}.

It remains to study the average order of $S_u(q;\v)$, 
as $\v$ ranges over an interval of length $P^\ve V$, centred upon a point $\v_0$.
Our investigation of this topic will adhere to the basic approach
in \cite{41}. Beginning with the multiplicativity property of the sums
$S_u(q;\v)$, we have the following result.

\begin{lem}\label{lem:S-mult}
We have 
$$
S_{u}(rs,\v)=S_{u\bar{r}^2}(s,\bar{r}\v)S_{u\bar{s}^2}(r,\bar{s}\v),
$$
provided that $r,s$ are coprime and $\bar{r}, \bar{s}$ are any integers such
that $r\bar{r}+s\bar{s}=1.$   
\end{lem}

\begin{proof}
Now it follows from \cite[Lemma 10]{41} that
$$
T(a,rs;\v)=T(a\bar{s},r;\bar{s}\v)T(a\bar{r},s;\bar{r}\v),
$$
where 
\begin{equation}
  \label{eq:TT}
T(a,q;\v):=\sum_{\y\bmod{q}}e_q(ag(\y)+\v.\y).
\end{equation}
As $\alpha$ (resp. $\beta$) ranges over integers coprime to $r$ modulo
$r$ (resp. coprime to $s$ modulo $s$), so $a=\al s\bar{s}+\beta r
\bar{r}$ ranges over a set of residues modulo $q$ that are coprime to
$q$. Hence it follows that
\begin{align*}
S_{u}(rs,\v)&=
\sum_{\colt{\alpha \bmod{r}}{\hcf(\alpha,r )=1}}
\sum_{\colt{\beta \bmod{s}}{\hcf(\beta,s )=1}}
e_{rs}\big(u(\bar{\al} s\bar{s}+\bar{\beta} r
\bar{r})\big)
T(\alpha \bar{s},r;\bar{s}\v)T(\beta \bar{r},s;\bar{r}\v)\\
&=S_{u\bar{r}^2}(s,\bar{r}\v)S_{u\bar{s}^2}(r,\bar{s}\v),
\end{align*}
as required for the statement of the lemma.
\end{proof}

Lemma \ref{lem:S-mult} allows us to concentrate on the value of 
$S_{u}(q,\v)$  at prime power values of $q$. 
We begin by considering the sums $S_{u}(p,\v)$ for primes $p$ such
that $g_0$ remains non-singular modulo $p$. 
It follows from \cite[Lemma 4]{cubshort}
that 
\begin{equation}
  \label{eq:SP1}
S_u(p;\v) \ll p^{(n+1)/2},
\end{equation}
provided that $p\nmid u$. 
For the case in which $p\mid u$ the situation is more complicated.
Let $g^*\in \Z[\x]$ denote the dual form to $g_0$. Since $g_0$ is
non-singular, so it follows that $g^*$ is absolutely irreducible and
has degree $D$, with $2\leq D \ll 1$.  It now follows from
\cite[Lemma 5]{cubshort} that 
\begin{equation}
  \label{eq:SP2}
S_0(p;\v) \ll p^{(n+1)/2}(p,g^*(\v))^{1/2}.
\end{equation}
When $p$ is a prime such that $g_0$ is singular modulo $p$, so that
$p$ divides the discriminant $\D_{g_0}$ of $g_0$, it will suffice to employ the
trivial bound
\begin{equation}
  \label{eq:SP3}
|S_u(p^j;\v)| \leq p^{j(n+1)} \leq \hcf(p,\D_{g_0})^\theta
\ll \|g_0\|^\theta \leq H^{\theta},
\quad (j=1,2).
\end{equation}
To handle moduli involving $p^2$ when $p\nmid \D_{g_0}$ we
employ \cite[Lemma~7]{41}, giving 
\begin{equation}
  \label{eq:SP4}
S_u(p^2;\v) \ll p^{n+2}.
\end{equation}

Recall the decomposition $q=bc^2 d=b_1b_2^2c^2d$ in  \eqref{eq:bcd},
with $b_1,b_2^2,c^2d$ pairwise coprime.
Drawing together \eqref{eq:SP1}, \eqref{eq:SP2} and \eqref{eq:SP3}, we
conclude from Lemma \ref{lem:S-mult} that 
$$
S_u(b_1;\v)\leq  H^{\theta}
A^{\omega(b_1)}b_1^{(n+1)/2}\hcf(b_1,u,g^*(\v))^{1/2},
$$
for some absolute constant $A\ll 1$. Similarly, it follows
from \eqref{eq:SP3} and \eqref{eq:SP4} that
$$
S_u(b_2;\v)\leq H^{\theta} A^{\omega(b_2)}b_2^{n+2}.
$$
We may now conclude that
$$
S_u(q;\v)\ll H^{\theta} b^{(n+1)/2+\ve}b_2\hcf(b_1,u,g^*(\v))^{1/2}
|S_{u\bar{b}^2}(c^2d;\bar{b}\v)|,
$$
for any $\ve>0$, where $\bar{b}\in\Z$ is such that
$b\bar{b}\equiv 1 \bmod{c^2d}$.

It is now time to introduce the summation over
$\ma{v}$ in \eqref{eq:S2}, for given $\v_0\ll HP$.  
Let us write $\mcal{S}_1$ for the overall contribution from the case
in which either $u\neq 0$ or else $g^*(\v)\neq 0$ in the summation
over $\v$. We write $\mcal{S}_2$ for the corresponding contribution
from the case $u=0$ and $g^*(\v)= 0$ in the summation over $\v$. 
It follows that
\begin{align*}
\mcal{S}_1 \ll H^{\theta} \max_{0< N \ll
  (HP)^{\theta}} 
q^{\ve} b^{(n+1)/2}\hcf(b_1,N)^{1/2}b_2
\sum_{|\v-\v_0|\leq P^\ve V} |S_{u\bar{b}^2}(c^2d;\bar{b}\v)|.
\end{align*}
Observing that 
$$
|S_{u\bar{b}^2}(c^2d;\bar{b}\v)| \leq \sum_{\colt{a
    \bmod{c^2d}}{\hcf(a,c^2d)=1}}|T(a,c^2d;\bar{b}\v)|, 
$$
where $T(a,c^2d;\bar{b}\v)$ is given by \eqref{eq:TT}, 
we may now combine \cite[Lemmas~11, 15 and 16]{41} in the manner indicated
at the close of \cite[\S 5]{41}, in order to conclude that 
\begin{equation}
  \label{eq:SS1}
\mcal{S}_1 \ll  H^{\theta}
q^{n/2+1+\ve}W^n
\max_{0< N \ll
  (HP)^{\theta}}\frac{\hcf(b_1,N)^{1/2}}{b_1^{1/2}}.
\end{equation}
Here $W$ is given by \eqref{eq:W}, and we have taken 
$
\max_{1\leq  i\leq n}r_i \ll  H^\theta
$ 
in \cite{41}.

In order to estimate $\mcal{S}_2$, in which case $u=0$, 
we begin as above with the observation that 
$$
\mcal{S}_2 \ll H^{\theta} 
q^{\ve} b^{n/2+1}\sum_{\colt{|\v-\v_0|\leq P^\ve V}{g^*(\v)=0}} |S_0(c^2d;\bar{b}\v)|.
$$
By dropping the condition $g^*(\v)=0$ we may certainly re-apply the work leading
to \eqref{eq:SS1}, giving 
\begin{equation}
  \label{eq:SS2a}
\mcal{S}_2 \ll  H^{\theta}
q^{n/2+1+\ve} W^n.
\end{equation}
Alternatively, we apply 
\cite[Lemma 11]{41} to deduce that 
\begin{align*}
\mcal{S}_2 
&\ll H^{\theta} 
q^{\ve}b^{n/2+1}(c^2d)^{n/2}
\sum_{\colt{a \bmod{c^2d}}{\hcf(a,c^2d)=1}}
\sum_{\colt{|\v-\v_0|\leq P^\ve V}{g^*(\v)=0}} 
\sum_{\colt{\ma{a}\bmod{c}}{
c\mid (a\nabla g(\ma{a})+\v)}}M_d(\ma{a})^{1/2}\\
&\ll H^{\theta} 
q^{n/2+1+\ve}\max_{\colt{a \bmod{c^2d}}{\hcf(a,c^2d)=1}}
\sum_{\colt{|\v-\v_0|\leq P^\ve V}{g^*(\v)=0}} 
\sum_{\colt{\ma{a}\bmod{c}}{
c\mid (a\nabla g(\ma{a})+\v)}}M_d(\ma{a})^{1/2},
\end{align*}
where
$$
M_d(\x):=
\#\big\{\y \bmod{d}: \nabla^2 g(\x)\y \equiv \ma{0} \mod{d}\big\}.
$$
We proceed to employ the following result to estimate the number of
available $\v$.

\begin{lem}\label{lem:g*}
We have
$$
\#\Big\{\ma{v}\in\Z^n: 
\begin{array}{l}
|\ma{v}-\ma{v}_0|\leq X, ~g^*(\ma{v})=0,\\
\ma{v}\equiv \ma{a} \bmod{q}
\end{array}
\Big\}\ll 1+ \Big(\frac{X}{q}\Big)^{n-3/2+\ve},
$$
uniformly in $\ma{v}_0, \ma{a}$ and the coefficients of $g^*$.
\end{lem}

\begin{proof}
Let $N(q;X)$ denote the quantity that is to be estimated. 
Dropping the condition that $g^*(\ma{v})=0$, we 
trivially have $N(q;X)=O(1)$ if $q>X$.  Assume henceforth that $q \leq X$, and
write $\v=\ma{a}+q\ma{w}$ for $\ma{w}\in\Z^n$. Since $|\v-\v_0|\leq
X$, so it follows that $|\w-\w_0|\ll X/q$, where
the components of $\w_0$ are obtained by taking the integer part of
the components of $q^{-1}(\v_0-\ma{a})$.
On writing $\ma{z}=\w-\w_0$, we therefore deduce that
$$
N(q;X)\leq \#\big\{\ma{z}\in\Z^n: \ma{z}\ll X/q, ~f(\ma{z})=0\big\},
$$
where $f(\ma{z})=g^*(\ma{a}+q\w_0+q\ma{z})$.  

In view of the fact
$g^*$ is an absolutely irreducible form of degree at least $2$,
so it follows that the shifted polynomial $f$ must be absolutely
irreducible with degree at least $2$. Our new polynomial $f$ need not be
homogeneous, and has coefficients that depend on $\ma{a}, \ma{v}_0$
and $q$. We now appeal to a very general uniform bound due to Pila
\cite[Theorem~A]{pila}, which implies that the affine hypersurface
$f=0$ contains  
$$
\ll Y^{n-2+1/\deg f+\ve}
$$ 
integer points of height 
at most $Y$, for any $Y\geq 1$.  Furthermore, the implied constant in
this estimate is uniform in the coefficients of $f$. Once inserted
into our bound for 
$N(q;X)$, this therefore completes the proof of the lemma.
\end{proof}

Employing Lemma \ref{lem:g*} in our bound for $\mcal{S}_2$ we deduce that
\begin{align*}
\mcal{S}_2
&\ll H^\theta q^{n/2+1+\ve} 
\max_{\colt{a\bmod{c^2d}}{\hcf(a,c^2d)=1}}
\sum_{\ma{a}\bmod{c}}M_d(\ma{a})^{1/2}
\sum_{\tstack{|\v-\v_0|\leq P^\ve V}{g^*(\v)=0}{
c\mid (a\nabla g(\ma{a})+\v)}}1\\
&\ll H^\theta q^{n/2+1} P^\ve \Big(1+\frac{V}{c}\Big)^{n-3/2}
U(q),
\end{align*}
where 
$$
U(q):=\sum_{\ma{a}\bmod{c}}
M_d(\ma{a})^{1/2}.
$$
At this point it is worth comparing our investigation with the corresponding
argument in \cite{cubshort} and \cite{hb-10}. There
a refinement of \cite[Lemma 11]{41} is used in order to obtain 
a better estimate for $S_0(c^2d;\bar{b}\v)$ when $g^*(\v)=0$. This
ultimately leads to a version of the
above bound for $\mcal{S}_2$ with the additional constraint that
$c\mid g(\ma{a})$ in the definition of $U(q)$. If we denote this
quantity by $U^*(q)$, then  in the 
setting of homogeneous polynomials Heath-Brown shows that $U^*(q)\ll_g
c^{n-1+\ve}d^{1/2}$, provided that $n\geq 10$. A slightly weaker
estimate is achieved in \cite{cubshort}, but which is valid for
cubic polynomials that are not necessarily homogeneous. It turns out
that neither of these estimates is readily extended to smaller values
of $n$, but fortunately we have found it sufficient to work with
$U(q)$ instead. We can then combine \cite[Lemma 14]{41} with Cauchy's inequality to
deduce that
\begin{align*}
U(q)
\ll c^{n/2}\Big(\sum_{{\ma{a}\bmod{c}}} M_d(\ma{a})\Big)^{1/2}
&\ll c^{n/2}\Big(\frac{c^{n}}{d^{n}}\sum_{{\ma{a}\bmod{d}}}
M_d(\ma{a})\Big)^{1/2}\\
&\ll H^\theta c^{n+\ve}.
\end{align*}
Once inserted into the preceding bound for $\mcal{S}_2$, we deduce that
\begin{equation}\label{eq:SS2b}
\mcal{S}_2
\ll H^\theta q^{n/2+1} P^\ve M_2,
\end{equation}
in the notation of Proposition \ref{pro:S}.

To obtain a final estimate for $\mcal{S}_2$, we simply drop the condition
that $g^*(\ma{v})=0$ in the summation over
$\ma{v}$. In this way \cite[Lemma
16]{41} easily leads us to the conclusion that
\begin{equation}
  \label{eq:tue2}
\mcal{S}_2\ll H^\theta q^{n/2+1} P^\ve M_3,
\end{equation}
in the notation of Proposition \ref{pro:S}.  

Drawing together \eqref{eq:SS1}, \eqref{eq:SS2a},
\eqref{eq:SS2b}, and \eqref{eq:tue2} 
in \eqref{eq:S2}, we therefore complete the proof of
Proposition \ref{pro:S}.

\section{Proof of Proposition \ref{main'}}\label{s:finale}

In this section we establish Proposition \ref{main'}. 
Let $w\in\WW_n$ and let $g\in \Z[\x]$ be a cubic polynomial with
$s(g_0)=-1$ and $n \geq 5$.
Taking \eqref{eq:Nw'} as our starting point, with the choice
$
Q=P^{3/2},
$
we have 
$$
N_{w}(g;P)=T_1+O(T_2),  
$$
where
$$
T_1:=\sum_{q\leq Q}\int_{\frac{-1}{qQ}}^{\frac{1}{qQ}} \SSS_0(q;z)\d z,
\quad
T_2:=\frac{1}{P^3}\sum_{q\leq Q}\sum_{|u|\leq
  \frac{q}{2}}\frac{\max_{\frac{1}{2}\leq qQ|z|\leq
    1}|\SSS_u(q;z)|}{1+|u|}.
$$
Let $P\geq 1$ and $H\geq \|g\|_P$. Throughout this section we may
assume that $P\geq H$, since the alternative hypothesis simply contributes
$O(H^n)$ to $N_{w}(g;P)$, which is satisfactory for Proposition \ref{main'}.

We will consider the contribution to $T_1,T_2$ from $q$ restricted to
lie in certain intervals. Write $q=b_1b_2^2c^2d$, where $b_1,b_2,d$ are given by
\eqref{eq:bcd}. Let $R,R_0,\ldots,R_3\geq 1/2$ and $t>0$. Then we will
write $\Sigma_{i}(R,\ma{R};t)$ for the overall contribution to
$T_i$, for $i=1,2$, from those $q,z$ for which 
\begin{equation}
  \label{eq:s0s1}
\begin{split}
  R<q\le 2R, \quad  
  R_0<b_1\le 2R_0, \quad  
  R_1<b_2\le 2R_1,\\
  R_2<c\le 2R_2 \quad
  R_3<d\le 2R_3,
\end{split}
\end{equation}
and 
$$
  t<|z|\leq 2t.
$$
Our plan will be to show that
\begin{equation}
  \label{eq:i}
\Sigma_{i}(R,\ma{R};t) \ll H^\theta P^{n-2+\ve},
\end{equation}
for $i=1,2$, under the assumption that $n \geq 5$ and
$s(g_0)=-1$. Once summed over $O((\log P)^6)$ dyadic intervals for 
$R,\ma{R}$ and $t$, this will clearly suffice to complete the
proof of Proposition \ref{main'}.

Recall that $d\mid c$. Thus $\Sigma_i(R,\ma{R};t)=0$ for $i=1,2$, unless 
\begin{equation}
  \label{eq:s0s1'}
  R_3\ll R_2, \quad R\ll R_0R_1^2R_2^2R_3\ll R \leq P^{3/2}.
\end{equation}
Similarly, it is clear that $\Sigma_i(R,\ma{R};t)=0$ unless
\begin{equation}
  \label{eq:dy-t}
  (RP^{3/2})^{-1} \geq t \gg \left\{
\begin{array}{ll}
0, & \mbox{if $i=1$},\\
(RP^{3/2})^{-1}, & \mbox{if $i=2$}.
\end{array}
\right.
\end{equation}
The following simple result will be useful in our work.

\begin{lem}\label{lem:count_s0s1}
We have
$$
\#\{q=b_1b_2^2c^2d:\,\mbox{\eqref{eq:s0s1} holds}\}
\ll R_0R_1R_2^{1/2}R_3^{1/2}.
$$
\end{lem}

\begin{proof}
It is clear that we have to count the number of quadruples
$(b_1,b_2,c,d)$ for which $d\mid c$ and \eqref{eq:s0s1} holds.  
The number of choices for $b_1$ and $b_2$ is $O(R_0R_1)$.  To
count the possible pairs $c,d$ recall from \eqref{eq:bcd} that there 
exist a positive integer $d_0$ such that
$d_0 \mid d$ and $d_0^{-1}d^{-1}c$ is square-full.
Hence, for fixed values of $d$, the number of 
available choices for $c$ is
$$
\ll \sum_{d_0 \mid d}\Big(\frac{R_2}{d_0R_3}\Big)^{1/2}.
$$ 
On summing over values of $d$, we deduce that the overall number of
choices for $c,d$ is
$$
\ll R_2^{1/2}R_3^{-1/2}\sum_{d\ll R_3}\sum_{d_0 \mid d} 
\frac{1}{d_0^{1/2}}
\ll R_2^{1/2}R_3^{1/2}.
$$ 
This suffices for the proof of Lemma \ref{lem:count_s0s1}.
\end{proof}

Our main tool in bounding $\Sigma_1(R,\ma{R};t)$ and $\Sigma_2(R,\ma{R};t)$ will be Proposition
\ref{pro:S}, but this will be supplemented with Proposition
\ref{pro:Sbirch} to handle certain awkward ranges of $R$.

\subsection{Estimating $\Sigma_2(R,\ma{R};t)$ }

We begin with our treatment of $\Sigma_2(R,\ma{R};t)$.
Since the size of $t$ is effectively determined by \eqref{eq:dy-t}, so
it will be convenient to write
$\Sigma_2(R,\ma{R})=\Sigma_2(R,\ma{R};t)$ throughout this section.

It follows from Proposition \ref{pro:S} that
\begin{align*}
\Sigma_2(R,\ma{R})
&\ll H^{\theta}P^{n-3+\ve}
\sum_{q}R^{1-n/2}\max_{|z|\asymp (RQ)^{-1}}(W^n M_1 + \min\{M_2,M_3\}),
\end{align*}
where $W, M_1,M_2,M_3$ are as in the statement of the proposition and the
summation over $q$ is over all $q=b_1b_2^2c^2d$ such that
$q,b_1,b_2,c, d$ are constrained to lie in the dyadic ranges \eqref{eq:s0s1}.
Let us write $\Sigma_{2,a}$ for the overall contribution to the right
hand side from the term involving $W^nM_1$, and 
$\Sigma_{2,b}$ for the corresponding contribution from the term 
involving $\min\{M_2,M_3\}$. Thus we have 
\begin{equation}
  \label{eq:65.1}
\Sigma_2(R,\ma{R})\ll  \Sigma_{2,a}+\Sigma_{2,b}.
\end{equation}

We begin by estimating $\Sigma_{2,a}$, for which it is convenient to
note that
\begin{equation}
  \label{eq:65.3}
V \asymp R^{1/2}P^{-1/4}, \quad 
W\ll  R^{1/2}P^{-1/4}+  (R_2^2R_3)^{1/3},
\end{equation}
for any $z$ such that $|z|\asymp (RQ)^{-1}$, with
$Q=P^{3/2}$. Now it is trivial to see that
$$
\sum_{b\leq B}\hcf(b,N)\ll \tau(N)B \ll N^{\ve}B,
$$
for any $N\in\N$ and $B\geq 1$.  By adjusting the proof of Lemma
\ref{lem:count_s0s1} slightly it therefore follows that
$$
\sum_q \hcf(b_1,N)^{1/2} \ll N^{\ve} R_0R_1R_2^{1/2}R_3^{1/2},
$$
for any $N\in\N$.  Bringing this all together we conclude that
\begin{align*}
\Sigma_{2,a}
&\ll H^{\theta}\frac{P^{n-3+\ve}}{R_0^{1/2}R^{n/2-1}} 
\max_{0<N\ll (HP)^{\theta}}
\sum_{q} \hcf(b_1,N)^{1/2}
\max_z W^n\\
&\ll H^{\theta}\frac{P^{n-3+\ve}R_0^{1/2}R_1R_2^{1/2}R_3^{1/2}}{R^{n/2-1}} 
\max_z W^n\\
&\ll H^{\theta}\frac{P^{n-3+\ve}}{R^{n/2-3/2}} 
(R^{1/2}P^{-1/4}+  (R_2^2R_3)^{1/3})^n,
\end{align*}
since $R_0\ll R/(R_1^2R_2^2R_3)$ by \eqref{eq:s0s1'}.
Our aim is to show that 
\begin{equation}
  \label{eq:65.2}
  \Sigma_{2,a}\ll H^{\theta} P^{n-2+\ve},
\end{equation}
provided that $n\geq 5$. We have two terms to consider in our estimate
for $\Sigma_{2,a}$. Beginning with the term 
involving $R^{1/2}P^{-1/4}$, we obtain the contribution 
$$
\ll H^{\theta}P^{3n/4-3+\ve}R^{3/2} \ll H^{\theta}P^{3n/4-3/4+\ve},
$$
since $R\leq P^{3/2}$. This is satisfactory for $n\geq 5$.  Finally, the term involving 
$(R_2^2R_3)^{1/3}$ contributes
\begin{align*}
&\ll
H^{\theta}\frac{P^{n-3+\ve}(R_2^2R_3)^{n/3}}{R^{n/2-3/2}}
\ll
H^{\theta}P^{n-3+\ve}R^{3/2-n/6}.
\end{align*}
since $R_2^2R_3\ll R$.  When $n\geq 9$ the exponent of $R$ is
non-positive, which clearly yields a satisfactory contribution. 
When $5\leq n\leq 8$, we obtain the contribution $O(H^{\theta} P^{3n/4-3/4+\ve})$
by taking $R\leq P^{3/2}$. This completes the proof of \eqref{eq:65.2}.

We now turn to the task of estimating $\Sigma_{2,b}$, for which we
want to show that
\begin{equation}
  \label{eq:65.2'}
  \Sigma_{2,b}\ll H^{\theta} P^{n-2+\ve},
\end{equation}
provided that $n\geq 5$.  Once combined with \eqref{eq:65.2} in
\eqref{eq:65.1} this will be enough to establish \eqref{eq:i} in the
case $i=2$. We will need to supplement our estimate with Proposition 
\ref{pro:Sbirch}. A little thought reveals that 
\begin{align*}
\Sigma_{2,b}
&\ll H^{\theta}P^{n-3+\ve}
R\sum_{q}   \min\Big\{ P^{-3n/16}   , R^{-n/2}\max_{z}\min\{M_2,M_3\}\Big\},
\end{align*}
where $M_2,M_3$ are as in the statement of Proposition
\ref{pro:S}, and the maximum is over $z$ such that $|z|\asymp
(RQ)^{-1}$. In particular \eqref{eq:65.3} holds
in the definitions of $M_2,M_3$. 

Suppose first that $V\geq R_2$. Then we take $\min\{M_2,M_3\}\leq
M_2$, in order to conclude from \eqref{eq:s0s1'} and Lemma
\ref{lem:count_s0s1} that
\begin{align*}
\Sigma_{2,b}
&\ll H^{\theta}P^{n-3+\ve}
R^{1-n/2}R_0R_1R_2^{n+1/2}R_3^{1/2}\Big(\frac{V}{R_2}\Big)^{n-3/2}\\ 
&\ll H^{\theta}P^{n-3+\ve}
R^{2-n/2}V^{n-3/2}\\ 
&\ll H^{\theta}P^{3n/4-21/8+\ve}R^{5/4}\\
&\ll H^{\theta}P^{3n/4-3/4+\ve}.
\end{align*}
This is satisfactory for $n\geq 5$. 

Suppose now that $(R_2^2R_3)^{1/3}\leq  V< R_2$. Then we may take 
$$
\min\{M_2,M_3\}\leq M_2^{3/10}M_3^{7/10} \ll 
R_2^{3n/10} \Big(\frac{R^{1/2}}{P^{1/4}}\Big)^{7n/10}
$$
in the above. Lemma \ref{lem:count_s0s1} and
\eqref{eq:s0s1'} together reveal that 
\begin{align*}
\Sigma_{2,b}
&\ll H^{\theta}P^{n-3+\ve}
\frac{R^{2-n/2}}{R_2^{3/2}R_3^{1/2}} \min\{M_2,M_3\}
= H^{\theta}P^{n-2+\ve}E_n,
\end{align*}
where
\begin{align*}
E_n&=P^{-1-7n/40}R^{2-3n/20}R_2^{3n/10-3/2}.
\end{align*}
We wish to show that $E_n\ll 1$ for $n\geq 5$. But clearly
$3n/10-3/2\geq 0$ for $n$ in this range, whence we may take $R_2\ll R^{1/2}$ in
this estimate. It follows that
$$
E_n\ll P^{-1-7n/40}R^{5/4} \ll P^{7/8-7n/40} \ll 1,
$$
for $n\geq 5$, as required. 

Turning to the case in which $V< (R_2^2R_3)^{1/3}$, we note that
$$
M_3\ll \Big(\frac{R_2^2R_3}{V}\Big)^{n/2} 
\ll \frac{P^{n/8}(R_2^2R_3)^{n/2}}{R^{n/4}}
$$ 
in the statement of Proposition \ref{pro:S}. Hence
\begin{align*}
\Sigma_{2,b}
&\ll H^{\theta}
\frac{P^{n-3+\ve}R^2}{R_2^{3/2}R_3^{1/2}} \min\Big\{ \frac{1}{P^{3n/16}}   , 
\frac{P^{n/8}(R_2^2R_3)^{n/2}}{R^{3n/4}}
\Big\}\\
&\ll H^{\theta}
P^{n-3+\ve} \min\Big\{ \frac{R^2}{P^{3n/16}}   , 
\frac{P^{n/8}}{R^{n/4-4/3}}
\Big\},
\end{align*}
since $R_2^{3/2}R_3^{1/2}\gg (R_2^2R_3)^{2/3}$.
When $n=5$, so that $n/4-4/3<0$, we take $R\leq P^{3/2}$ to deduce
that
\begin{align*}
\Sigma_{2,b}
&\ll H^{\theta}
P^{21/8+\ve} (P^{3/2})^{1/12} \ll H^{\theta}
P^{11/4+\ve}\ll H^{\theta} P^{3},
\end{align*}
which is satisfactory. When $n\geq 6$ we apply
the bound coming from Weyl differencing when $R<P$, and the
bound coming from Poisson summation when $R\geq
P$. This yields
\begin{align*}
\Sigma_{2,b}
&\ll H^{\theta}P^\ve \big( P^{13n/16-1}+P^{7n/8-5/3}\big),
\end{align*}
which is satisfactory for $n\geq 6$. This completes the proof of
\eqref{eq:65.2'}.

\subsection{Estimating $\Sigma_1(R,\ma{R};t)$ }

It follows from Proposition \ref{pro:S} and the argument in Lemma \ref{lem:count_s0s1}
that
\begin{align*}
\Sigma_1(R,\ma{R};t)
&\ll H^{\theta}P^{n+\ve}t \Big(\frac{R^{3/2-n/2} W^n}{R_2^{1/2}} + 
\frac{R^{2-n/2} }{R_1R_2^{3/2}R_3^{1/2}}
\min\{M_2,M_3\}\Big),
\end{align*}
where $W, M_2,M_3$ are as in the statement of the proposition, but
with individual variables replaced by appropriate lower or upper
bounds corresponding to the interval that the variable is assumed to
lie in. Let us write $\Sigma_{1,a}$ for the overall contribution to the right
hand side from the first term, and 
$\Sigma_{1,b}$ for the corresponding contribution from the second
term. In order to establish \eqref{eq:i} with $i=1$, it will 
suffice to show that 
\begin{equation}
  \label{eq:56.1}
\max\{\Sigma_{1,a}, \Sigma_{1,b} \}\ll H^\theta P^{n-2+\ve},
\end{equation}
for $n\geq 5$.

Let us begin by estimating $\Sigma_{1,a}$, for which we have
\begin{equation}
  \label{eq:vv}
V\asymp \left\{
\begin{array}{ll}
R/P, &\mbox{if $t<P^{-3}$},\\
Rt^{1/2}P^{1/2}, & \mbox{if $t\geq P^{-3}$},
\end{array}
\right.
\end{equation}
and $W\ll   V +(R_2^2R_3)^{1/3}.$
When $t\geq P^{-3}$ the term involving $V$  makes the contribution
\begin{align*}
\ll H^{\theta}P^{n+\ve}t R^{3/2-n/2} (Rt^{1/2}P^{1/2})^n 
&\ll H^{\theta}P^{3n/2+\ve}t^{1+n/2} R^{3/2+n/2}\\
&\ll H^{\theta}P^{3n/4-3/4+\ve}
\end{align*}
to $\Sigma_{1,a}$, since $t\leq (RP^{3/2})^{-1}$. 
This is satisfactory for $n\geq 5$.
Likewise, when $t<P^{-3}$, one obtains a satisfactory contribution.
To handle the contribution from the term involving 
$(R_2^2R_3)^{1/3}$ we will need to supplement our estimate with an
application of Proposition \ref{pro:Sbirch}, in addition to differentiating according to the size
of $t$. Suppose first that $t\geq P^{-3}$. Then we have the overall contribution
\begin{align*}
&\ll H^{\theta}P^{n+\ve} 
 \min\Big\{\frac{R^{3/2-n/2} t (R_2^2R_3)^{n/3}}{R_2^{1/2}},
\frac{R^{2-n/8}t^{1-n/8}}{R_2^{3/2}R_3^{1/2}P^{3n/8}}\Big\}\\
&\ll H^{\theta}P^{n+\ve} 
 \min\Big\{R^{3/2-n/2} t (R_2^2R_3)^{n/3-1/6},
\frac{R^{2-n/8}t^{1-n/8}}{(R_2^2R_3)^{2/3}P^{3n/8}}\Big\},
\end{align*}
since $R_3\ll R_2$. We apply the basic inequality $\min\{A,B\}\leq
A^{1/3}B^{2/3}$ to derive the overall contribution
$O(H^{\theta}P^{n-2+\ve}E_n)$, with
$$
E_n=P^{2-n/4}t^{1-n/12}R^{11/6-n/4}(R_2^2R_3)^{n/9-1/2}.
$$
Suppose first that $5\leq n\leq 12$. Then we may take $t\leq
(RP^{3/2})^{-1}$ to deduce that
\begin{align*}
E_n\leq 
P^{1/2-n/8}R^{5/6-n/6}(R_2^2R_3)^{n/9-1/2}\ll 
P^{1/2-n/8}R^{1/3-n/18},
\end{align*}
whence $E_n\ll 1$. Alternatively, when $n\geq 13$ we have
\begin{align*}
E_n\leq 
P^{-1}R^{11/6-n/4}(R_2^2R_3)^{n/9-1/2}\ll P^{-1}R^{4/3-5n/36}\ll 1,
\end{align*}
since $t\geq P^{-3}.$

So far we have established a satisfactory bound for $\Sigma_{1,a}$ 
under the assumption that $t\geq P^{-3}$. When $t<P^{-3}$, we easily obtain the overall
contribution 
\begin{align*}
\ll H^{\theta}P^{n+\ve} 
 \frac{R^{3/2-n/2} t (R_2^2R_3)^{n/3}}{R_2^{1/2}} &\ll
 H^{\theta}P^{n-3+\ve} 
R^{3/2-n/2}  (R_2^2R_3)^{n/3-1/6} \\
&\ll
 H^{\theta}P^{n-3+\ve} 
R^{4/3-n/6}.
\end{align*}
The exponent of $R$ is non-positive when $n\geq 8$, in which case the bound
is clearly satisfactory. When $5\leq n\leq 7$, we take $R\leq P^{3/2}$
to obtain the satisfactory contribution $O(H^\theta P^{3n/4-1+\ve})$.
This establishes the bound for $\Sigma_{1,a}$ recorded in
\eqref{eq:56.1}, for $n\geq 5$. 

We now turn to the task of estimating $\Sigma_{1,b}$. We have 
\begin{equation}\label{eq:alg}
\Sigma_{1,b}
\ll H^{\theta}P^{n+\ve}t 
\frac{R^{2-n/2} }{R_2^{3/2}R_3^{1/2}}
\min\{M_2,M_3,M_4\},
\end{equation}
with
\begin{equation}\label{eq:alg'}
M_2=R_2^n\Big(1+\frac{V}{R_2}\Big)^{n-3/2}, \quad 
M_3=V^n\Big(1+\frac{R_2^2R_3}{V^3}\Big)^{n/2},
\end{equation}
and 
\begin{equation}\label{eq:alg''}
M_4=R^{3n/8}\min\big\{1, (tP^3)^{-n/8}\big\}.
\end{equation}
Here $M_4$ arises from an application of Proposition \ref{pro:Sbirch}
and $V$ satisfies \eqref{eq:vv}. 
Let us begin by handling the
case in which $t\geq P^{-3}$, so that $V \asymp Rt^{1/2}P^{1/2}.$  
Suppose first that $V\geq R_2$. Then we take $\min\{M_2,M_3,M_4\}\leq
M_2$, in order to conclude that
\begin{align*}
\Sigma_{1,b}
&\ll H^{\theta}P^{n+\ve}t 
R^{2-n/2} R_2^{n-3/2}\Big(\frac{Rt^{1/2}P^{1/2}}{R_2}\Big)^{n-3/2}\\
&\ll H^{\theta}P^{3n/2-3/4+\ve}t^{n/2+1/4} 
R^{n/2+1/2}\\
&\ll H^{\theta}P^{3n/4-9/8+\ve}R^{1/4}\\
&\ll H^{\theta}P^{3n/4-3/4+\ve}.
\end{align*}
This is satisfactory for $n\geq 5$. 

Suppose now that $(R_2^2R_3)^{1/3}\leq  V< R_2$. We take 
$\min\{M_2,M_3,M_4\}\leq M_2^{3/10}M_3^{7/10}$ in \eqref{eq:alg}. This gives
$\Sigma_{1,b} \ll H^{\theta}P^{n-2+\ve}E_n$, where 
$$
E_n\ll P^{1/2-7n/20}R^{1/4} \ll P^{7/8-7n/20}\ll 1
$$ 
for $n\geq 5$, since $t\leq (RP^{3/2})^{-1}$.
Finally we consider the case $V< (R_2^2R_3)^{1/3}$. 
In this setting we have $M_2\ll R_2^n$ and $M_3\ll (R_2^2R_3/V)^{n/2}$
in \eqref{eq:alg'}, and $M_4\leq R^{3n/8}(tP^3)^{-n/8}$ in \eqref{eq:alg''}.
Taking $\min \{A,B,C\}\leq A^{1/10}B^{1/5}C^{7/10}$ in \eqref{eq:alg},
we therefore deduce that $\Sigma_{1,b} \ll H^{\theta}P^{n-2+\ve}E_n$, with
\begin{align*}
E_n=P^{2-5n/16}t^{1-11n/80}R^{2-27n/80}R_2^{3n/10-3/2}R_3^{n/10-1/2}.
\end{align*}
Suppose first that $n\leq 7$, so that $1-11n/80\geq 0$. Then the upper
bound $t\leq (RP^{3/2})^{-1}$ gives
\begin{align*}
E_n 
&\leq P^{1/2-17n/160}R^{1-n/5}R_2^{3n/10-3/2}R_3^{n/10-1/2}\\
&\ll P^{1/2-17n/160}R^{1/4-n/20}\\
&\ll 1,
\end{align*}
since $R_2^{3n/10-3/2}R_3^{n/10-1/2} \ll R^{3n/20-3/4}$ for $n\geq
5$.  When $n\geq 8$ we instead take $t\geq P^{-3}$ in the above, obtaining
$
E_n\leq P^{-1+n/10}R^{5/4-3n/16}\ll 1,
$
when $n\leq 10.$ Finally, when $n\geq 11$, we instead take
$\min\{A,B,C\}\leq C$ in the above to get a satisfactory contribution.

In order to complete the proof of \eqref{eq:56.1} it remains 
to show that $\Sigma_{1,b}\ll H^{\theta} P^{n-2+\ve}$ when $n \geq 5$
and $t<P^{-3}$.
In particular we have $V\asymp R/P$, by \eqref{eq:vv}, and it now
follows from \eqref{eq:alg} that
\begin{align*}
\Sigma_{1,b}
\ll& H^{\theta}P^{n-3+\ve}
\frac{R^{2-n/2}}{R_2^{3/2}R_3^{1/2}}
\min\{M_2,M_3,M_4\},
\end{align*}
with $M_2,M_3,M_4$ being given by \eqref{eq:alg'} and \eqref{eq:alg''}.
When $V\geq R_2$, we take $M_2$ in the minimum, giving 
$$
\Sigma_{1,b}
\ll H^{\theta}P^{n-3+\ve} R^{n/2+1/2}/P^{n-3/2} \ll H^\theta P^{3n/4-3/4+\ve}.
$$
This is satisfactory for $n \geq 5$. 
When $(R_2^2R_3)^{1/3}\leq V <R_2$ we take 
$\min\{M_2,M_3,M_4\}\leq M_2^{3/10}M_3^{7/10} $ to deduce that
$\Sigma_{1,b} \ll H^{\theta}P^{n-2+\ve}E_n$, with 
$$
E_n \ll \frac{R^{n/5+2}R_2^{3n/10-3/2}}{P^{7n/10+1} } \ll
\frac{R^{7n/20+5/4}}{P^{7n/10+1}} \ll P^{7/8-7n/40}\ll 1,
$$
for $n \geq 5$.

Finally we must deal with the case $V<(R_2^2R_3)^{1/3}$, in which
setting 
\begin{align*}
\Sigma_{1,b}
&\ll H^{\theta}
\frac{P^{n-3+\ve}R^{2-n/2}}{R_2^{3/2}R_3^{1/2}} \min\Big\{ R_2^{n},
\Big(\frac{P R_2^2R_3}{R}\Big)^{n/2}, R^{3n/8}\Big\}.
\end{align*}
We use the inequality  $\min \{A,B,C\}\leq A^{1/10}B^{11/30}C^{8/15}$
to deduce that $\Sigma_{1,b} \ll H^{\theta}P^{n-2+\ve}E_n$, with
\begin{align*}
E_n&=P^{11n/60-1}R^{2-29n/60}R_2^{7n/15-3/2}R_3^{11n/60-1/2}\\
&\ll P^{11n/60-1}R^{3/2-3n/10}R_2^{n/10-1/2}.
\end{align*}
In particular we have $E_5\ll P^{-1/12}\ll 1$. When $n\geq 6$ we have
$$
E_n\ll P^{11n/60-1}R^{5/4-n/4}.
$$
This is clearly $O(1)$ when $R\geq P^{7/10}$ and $6\leq
n\leq 15$. Assume now that $n\geq 16$, or else $6\leq
n\leq 15$ and $R< P^{7/10}$. Then we take 
$\min\{A,B,C\}\leq C$ in the above estimate instead, obtaining
$\Sigma_{1,b} \ll H^{\theta}P^{n-2+\ve}E_n$, but this time with 
$$
E_n=P^{-1}R^{2-n/8}.
$$
If $6\leq n\leq 15$ and $R< P^{7/10}$ then clearly $E_n\leq
P^{2/5-7n/80}\ll 1$. Alternatively, if $n\geq 16$ then $E_n\ll
P^{-1}\ll 1$.  This completes the proof of Proposition \ref{main'}.

\end{document}